\documentclass[11pt]{article}

\usepackage{amsmath,amssymb}
\usepackage{theorem}

\makeatletter
\@addtoreset{equation}{section}
\makeatother

\newtheorem{theorem}{Theorem}[section]
\newtheorem{proposition}[theorem]{Proposition}
\newtheorem{lemma}[theorem]{Lemma}

\newtheorem{definition}[theorem]{Definition}
{\theorembodyfont{\rmfamily}

\newtheorem{remark}[theorem]{Remark}}

\title{\Huge\bf Generic symmetries of the Laurent extension of quantum plane}

\author{\Large Sergey Sinel'shchikov}
\date{\small\it Mathematics Division, B. Verkin
Institute for Low Temperature Physics and Engineering,
\\ National Academy of Sciences of Ukraine
\\ 47 Lenin Ave., 61103 Kharkov, Ukraine
\\ \rm E-mail: sinelshchikov@ilt.kharkov.ua}

\begin{document}
\large

\maketitle

\begin{quotation}\small
A list of generic $U_q(\mathfrak{sl}_2)$-module algebra structures on the Laurent polynomial algebra over the quantum plane with uncountably many isomorphism classes is produced. Also, a complete list of such structures is presented in which the action of Cartan generator of $U_q(\mathfrak{sl}_2)$ does not reduce to multiplying $x$ and $y$ (the generators of quantum plane) by constants.

\vspace{3ex} {\it Key words}: quantum universal enveloping algebra, Hopf
algebra, Laurent polynomial, weight.

{\it Mathematics Subject Classification 2010}: 81R50, 17B37.
\end{quotation}

\maketitle

\bigskip

\section{Introduction}
It is well known that the quantum plane $\mathbb{C}_q[x,y]$ admits a structure of $U_q(\mathfrak{sl}_2)$-module algebra (see, e.g., \cite{kassel}). In fact, it was a single selected structure which was implicit in various related topics and applications. The question on to what extent this structure is unique was initially raised in \cite{DS}. It has been established therein that there exists an uncountable family of non-isomorphic structures of $U_q(\mathfrak{sl}_2)$-module algebra on quantum plane, and a complete classification of those has been presented. The next step has been done in a work by S. Duplij, Y. Hong, and F. Li \cite{DHL}, where the structures of $U_q(\mathfrak{sl}_m)$-module algebra on a generalized quantum plane, a polynomial algebra in $n$ quasi-commuting variables, $m,n>2$, are considered. In all the above cases, once $m,n$ are fixed, the structures in question (to be abbreviated as symmetries) belong to finitely many series. Every such series is labelled by a pair (in the simplest case \cite{DS}; in the more general context of \cite{DHL} their quantity is $(m-1)n$) of the so called weight constants, which determine the action of Cartan generators of the quantum universal enveloping algebra on the generators of (generalized) quantum plane.

In our opinion, somewhat different generalization of the results of \cite{DS} as compared to \cite{DHL}, makes a separate interest. Namely, instead of increasing the number $m$ and $n$ of generators, we suggest to retain $m=n=2$, but to add the inverse elements $x^{-1}$ and $y^{-1}$ for the generators of (the standard) quantum plane. This way we obtain the Laurent polynomial algebra $\mathbb{C}_q[x^{\pm 1},y^{\pm 1}]$ over the quantum plane. Our research demonstrates that this newly formed quantum algebra constitutes a much more symmetric object than the standard quantum plane. More precisely, the symmetries listed in \cite{DS} could be produced by separating out those symmetries on the extended algebra $\mathbb{C}_q[x^{\pm 1},y^{\pm 1}]$, which leave invariant the subalgebra $\mathbb{C}_q[x,y]$. Our list of symmetries looks more regular than that of \cite{DS}.

Our approach anticipates a passage through some additional difficulties. The latter are related to the fact that the action of the Cartan generator on the extended algebra $\mathbb{C}_q[x^{\pm 1},y^{\pm 1}]$ we consider in this work, does not reduce in general to multiplying $x$ and $y$ by weight constants, as it was the case in \cite{DS}, \cite{DHL}. This new context, after introducing some Preliminaries, is a subject of Section \ref{nw}. The complete list of symmetries in which monomials are not weight vectors, is given by Theorem \ref{sigma=-I}.

 On the contrary to the latter rather poor collection of symmetries, we use Section \ref{ws} to present the so called generic symmetries in which all the monomials are weight vectors (see Theorem \ref{gener}). The collection of generic series is abundant in the sense that it splits into uncountably many isomorphism classes of symmetries. On the contrary to \cite{DS}, \cite{DHL}, the collection of pairs of weight constants involved in generic symmetries, is also uncountable. These pairs of weight constants, coming from the generic symmetries, contain all but a countable family of weight constants (some rational powers of $q$) which appear as weight constants for (non-generic) symmetries. The latter symmetries are going to be a subject of a subsequent work.

\bigskip

\section{Preliminaries}

Let $H$ be a Hopf algebra whose comultiplication is $\Delta$, counit is
$\varepsilon$, and antipode is $S$ \cite{abe}. Also let $A$ be a unital
algebra whose unit is $\mathbf{1}$. We will also use the Sweedler notation
$\Delta(h)=\sum_ih_i'\otimes h_i''$ \cite{sweedler}.

\begin{definition}\label{symdef}
By a structure of $H$-module algebra on $A$ we mean a homomorphism of
algebras $\pi:H\to\operatorname{End}_\mathbb{C}A$ such that

(i) $\pi(h)(ab)=\sum_i\pi(h_i')(a)\cdot\pi(h_i'')(b)$ for all $h\in H$,
$a,b\in A$;

(ii) $\pi(h)(\mathbf{1})=\varepsilon(h)\mathbf{1}$ for all $h\in H$.

The structures $\pi_1,\pi_2$ are said to be isomorphic if there exists an
automorphism $\Psi$ of the algebra $A$ such that
$\Psi\pi_1(h)\Psi^{-1}=\pi_2(h)$ for all $h\in H$.
\end{definition}

Throughout the paper we assume that $q\in\mathbb{C}\setminus\{0\}$ is not a root of $1$ ($q^n\ne 1$ for all non-zero integers $n$). Consider the
quantum plane which is a unital algebra $\mathbb{C}_q[x,y]$ with two
generators $x,y$ and a single relation
\begin{equation}\label{qpr}
yx=qxy.
\end{equation}

Let us complete the list of generators with two more elements $x^{-1}$,
$y^{-1}$, and the list of relations with
\begin{equation}\label{Lpr}
xx^{-1}=x^{-1}x=yy^{-1}=y^{-1}y=\mathbf{1}.
\end{equation}
The extended unital algebra $\mathbb{C}_q[x^{\pm 1},y^{\pm 1}]$ defined
this way is called the Laurent extension of quantum plane (more precisely,
the algebra of Laurent polynomials over quantum plane).

Given an integral matrix $\sigma=\begin{pmatrix}k & l\\ m &
n\end{pmatrix}\in SL(2,\mathbb{Z})$ and a pair of non-zero complex numbers
$(\alpha,\beta)\in(\mathbb{C}^*)^2$, we associate an automorphism
$\varphi_{\sigma,\alpha,\beta}$ of $\mathbb{C}_q[x^{\pm 1},y^{\pm 1}]$
determined on the generators $x$ and $y$ by
\begin{equation}\label{Auto}
\varphi_{\sigma,\alpha,\beta}(x)=\alpha x^ky^m;\qquad
\varphi_{\sigma,\alpha,\beta}(y)=\beta x^ly^n.
\end{equation}
A well-known result claims that every automorphism of $\mathbb{C}_q[x^{\pm
1},y^{\pm 1}]$ has the form \eqref{Auto}, and the group
$\operatorname{Aut}(\mathbb{C}_q[x^{\pm 1},y^{\pm 1}])$ of automorphisms of
$\mathbb{C}_q[x^{\pm 1},y^{\pm 1}]$ is just the semidirect product of its
subgroups $SL(2,\mathbb{Z})$ and $(\mathbb{C}^*)^2$ determined by setting
\begin{equation}\label{sp}
\sigma(\alpha,\beta)\sigma^{-1}=(\alpha,\beta)^\sigma
\stackrel{\operatorname{def}}{=}(\alpha^k\beta^m,\alpha^l\beta^n).
\end{equation}
\cite{KPS} (see also \cite{AD}, \cite{PLCCN}).

The quantum universal enveloping algebra $U_q\left(\mathfrak{sl}_2\right)$
is a unital associative algebra defined by its (Chevalley) generators
$\mathsf{k}$, $\mathsf{k}^{-1}$, $\mathsf{e}$, $\mathsf{f}$, and the
relations
\begin{align}
\mathsf{k}^{-1}\mathsf{k} &=\mathbf{1},\quad\mathsf{kk}^{-1}=\mathbf{1},&
\label{kk1}
\\ \mathsf{ke} &=q^2\mathsf{ek},&\label{ke}
\\ \mathsf{kf} &=q^{-2}\mathsf{fk},&\label{kf}
\\ \mathsf{ef}-\mathsf{fe} &=\dfrac{\mathsf{k}-\mathsf{k}^{-1}}{q-q^{-1}}.& \label{effe}
\end{align}

The standard Hopf algebra structure on $U_q(\mathfrak{sl}_2)$ is determined
by
\begin{align}
\Delta(\mathsf{k}) &=\mathsf{k}\otimes\mathsf{k}, & \label{k0}
\\ \Delta(\mathsf{e}) &=\mathbf{1}\otimes\mathsf{e}+
\mathsf{e}\otimes\mathsf{k}, & \label{def}
\\ \Delta(\mathsf{f}) &=\mathsf{f}\otimes\mathbf{1}+
\mathsf{k}^{-1}\otimes\mathsf{f}, & \label{def1}
\\ \mathsf{S}(\mathsf{k}) &=\mathsf{k}^{-1}, &
\mathsf{S}(\mathsf{e}) &=-\mathsf{ek}^{-1}, &
\mathsf{S}(\mathsf{f}) &=-\mathsf{kf}, &\label{S}
\\ \boldsymbol{\varepsilon}(\mathsf{k}) &=\mathbf{1}, &
\boldsymbol{\varepsilon}(\mathsf{e})
&=\boldsymbol{\varepsilon}(\mathsf{f})=0. &\label{eps}
\end{align}

\bigskip

\section{The symmetries with non-trivial \boldmath$\sigma$}\label{nw}

It should be observed that, given a $U_q(\mathfrak{sl}_2)$-module algebra
structure on $\mathbb{C}_q[x^{\pm 1},y^{\pm 1}]$ (to be referred to as a
symmetry or merely an action for brevity), the generator $\mathsf{k}$ acts
via an automorphism of $\mathbb{C}_q[x^{\pm 1},y^{\pm 1}]$, as one can
readily deduce from invertibility of $\mathsf{k}$, Definition
\ref{symdef}\hspace{-.35em}(i) and \eqref{k0}. In particular, every
symmetry determines uniquely a matrix $\sigma\in SL(2,\mathbb{Z})$ as in
\eqref{Auto}.

\medskip

\begin{remark}
It turns out that there exists a one-to-one correspondence between the
$U_q(\mathfrak{sl}_2)$-symmetries that leave invariant the subalgebra
$\mathbb{C}_q[x,y]$ and the $U_q(\mathfrak{sl}_2)$-symmetries on
$\mathbb{C}_q[x,y]$. One can readily restrict such symmetry of
$\mathbb{C}_q[x^{\pm 1},y^{\pm 1}]$ to $\mathbb{C}_q[x,y]$.

On the other hand, suppose we are given an arbitrary symmetry $\pi$ on
$\mathbb{C}_q[x^{\pm 1},y^{\pm 1}]$ (not necessarily leaving invariant
$\mathbb{C}_q[x,y]$). One has the following relations:
\begin{flalign}
\pi(\mathsf{k})(x^{-1}) \!&=\!(\pi(\mathsf{k})x)^{-1} &
\pi(\mathsf{k})(y^{-1}) \!&=\!(\pi(\mathsf{k})y)^{-1} & \label{ext_k}
\\ \pi(\mathsf{e})(x^{-1}) \!&=\!-x^{-1}(\pi(\mathsf{e})x)(\pi(\mathsf{k})x)^{-1}
\!\!\! & \pi(\mathsf{e})(y^{-1})
\!&=\!-y^{-1}(\pi(\mathsf{e})y)(\pi(\mathsf{k})y)^{-1}\!\!\! &
\label{ext_e}
\\ \pi(\mathsf{f})(x^{-1})
\!&=\!-(\pi(\mathsf{k}^{-1})x)^{-1}(\pi(\mathsf{f})x)x^{-1}\!\! &
\pi(\mathsf{f})(y^{-1})
\!&=\!-(\pi(\mathsf{k}^{-1})y)^{-1}(\pi(\mathsf{f})y)y^{-1}\!\!\!\!\!\!\!\!\!
& \label{ext_f}
\end{flalign}
Here \eqref{ext_k} is straightforward since $\pi(\mathsf{k})$ is an
automorphism; \eqref{ext_e} and \eqref{ext_f} are derivable by
`differentiating' (i.e., applying $\mathsf{e}$ and $\mathsf{f}$,
respectively, to) \eqref{Lpr}. Certainly, these relations remain true when $x$ or $y$ is replaced by an arbitrary invertible element.

Thus, given a symmetry on $\mathbb{C}_q[x,y]$, the relations \eqref{ext_k} -- \eqref{ext_f} determine a well-defined extension of it to the additional generators $x^{-1}$, $y^{-1}$, hence to $\mathbb{C}_q[x^{\pm 1},y^{\pm 1}]$.
\end{remark}

It might look like there should be a sharp difference with the picture discovered in \cite{DS}, due to the (conjectured) abundance of non-weight actions related to non-trivial matrices $\sigma$. However, it turns out that there exists only a small collection of such symmetries. Let us start with describing those actions.

Suppose we are given a matrix $\sigma=\begin{pmatrix}k & l\\ m &
n\end{pmatrix}\in SL(2,\mathbb{Z})$, and let $\lambda$ and $\mu$ be its
eigenvalues. We consider the two series of assumptions on $\sigma$:

(i) $\overline{\lambda}=\mu$ and (since $\lambda\mu=1$)
$|\lambda|=|\mu|=1$,

(ii) $\lambda,\mu\in\mathbb{R}$ and $\lambda,\mu\notin\{-1,1\}$.

One observes first that $\lambda+\mu\in\mathbb{Z}$. In the case (i),  this, together with other restrictions of (i), implies that the only possible values for $\lambda+\mu=\operatorname{tr}\sigma$ could be $0$, $\pm 1$, $\pm 2$. Thus, being intended to find the symmetries corresponding to the assumptions of (i), we need to consider separately the five subcases.

\medskip

(i-1) Suppose that $\operatorname{tr}\sigma=2$. This means that
$\lambda=\mu=1$.

The case when $\sigma$ is just the identity matrix $I$ is postponed till the next Section.

Let us consider the case of (matrix conjugate to a) Jordan block, that is, the eigenspace of $\sigma$ is one dimensional. We need the following

\begin{lemma}\label{Czsym}
The complete list of $U_q(\mathfrak{sl}_2)$-symmetries on the Laurent polynomial algebra $\mathbb{C}[z^{\pm 1}]$ of a single variable $z$ is as follows.

1). Let $\gamma\in\mathbb{C}\setminus\{0\}$ be such that $\gamma^{r-1}=q^2$ for some $r\in\mathbb{Z}$. There exists a one-parameter ($a\in\mathbb{C}\setminus\{0\}$) family of $U_q(\mathfrak{sl}_2)$-symmetries on $\mathbb{C}[z^{\pm 1}]$ given by
$$
\pi(\mathsf{k})(z)=\gamma z;\qquad\pi(\mathsf{e})(z)=\frac{a}{q^2-1}z^r;\qquad
\pi(\mathsf{f})(z)=q^3(\gamma-1)a^{-1}z^{2-r}.
$$
Additionally, there exist two more symmetries
$$
\pi(\mathsf{k})(z)=\pm z;\qquad\pi(\mathsf{e})(z)=\pi(\mathsf{f})(z)=0.
$$
All the symmetries with fixed $\gamma$ are isomorphic, e.g., to that with $a=1$. There exists an isomorphism between the symmetries with $\gamma$ and $\gamma^{-1}$. In all other cases the symmetries are non-isomorphic.

2). Let $\gamma\in\mathbb{C}\setminus\{0\}$. There exists a $U_q(\mathfrak{sl}_2)$-symmetry on $\mathbb{C}[z^{\pm 1}]$ given by
$$
\pi(\mathsf{k})(z)=\gamma z^{-1};\qquad\pi(\mathsf{e})(z)=\pi(\mathsf{f})(z)=0.
$$
All these symmetries are isomorphic, e.g., to that with $\gamma=1$.

The symmetries from 1) are non-isomorphic to those from 2).
\end{lemma}

{\bf Proof.} Since, with $\pi$ being a symmetry, $\pi(\mathsf{k})$ is an automorphism of $\mathbb{C}[z^{\pm 1}]$, and any automorphism of $\mathbb{C}[z^{\pm 1}]$ is given by either $z\mapsto\gamma z$ or $z\mapsto\gamma z^{-1}$, $\gamma\in\mathbb{C}\setminus\{0\}$ \cite{PLCCN}, we need to consider the two cases.

1). Let $\pi(\mathsf{k})(z)=\gamma z$.

Assume first $\pi(\mathsf{e})(z)\ne 0$. As a consequence of \eqref{ke} we have $\pi(\mathsf{ke})(z)=q^2\gamma\pi(\mathsf{e})(z)$, and with $\pi(\mathsf{e})(z)=\sum_ia_iz^i$, the assumption $a_r\ne 0$ implies $a_r\gamma^rz^r=q^2\gamma_rz^r$, hence
\begin{equation}\label{rd}
\gamma^{r-1}=q^2.
\end{equation}
Since $q$ is not a root of $1$, such $r\in\mathbb{Z}$ is unique.

So we establish that $\pi(\mathsf{e})(z)=az^r$, $a\in\mathbb{C}\setminus\{0\}$.

A completely similar argument allows one to deduce that $\pi(\mathsf{f})(z)=bz^{2-r}$. Here $b\in\mathbb{C}\setminus\{0\}$, because otherwise $\pi(\mathsf{f})$ is identically zero on $\mathbb{C}[z^{\pm 1}]$. With the latter assumption, we observe that \eqref{effe}, being applied to $z$, fails, as its l.h.s. vanishes, while its r.h.s. is non-zero, since $\gamma$, due to \eqref{rd}, is not a root of $1$ together with $q$. This contradiction demonstrates that $\pi(\mathsf{f})(z)\ne 0$.

It remains to use our formulas for $\pi(\mathsf{e})(z)$ and $\pi(\mathsf{f})(z)$ in applying \eqref{effe} to $z$, in order to compute the relation between $a$ and $b$. This requires two additional formulas
\begin{alignat}{2}
\pi(\mathsf{e})(z^p) &=\frac{\gamma^p-1}{\gamma-1}\,az^{p+r-1} &&=
\frac{\gamma^p-1}{\gamma-1}z^{p-1}\pi(\mathsf{e})(z),\label{e(z^p)}
\\ \pi(\mathsf{f})(z^p) &=\frac{\gamma^{-p}-1}{\gamma-1}\,bz^{p+1-r} &&=
\frac{\gamma^{-p}-1}{\gamma-1}z^{p-1}\pi(\mathsf{f})(z),\label{f(z^p)}
\end{alignat}
$p\in\mathbb{Z}$, whose proof is completely routine. This way we deduce that $ab=\frac{q^2(\gamma-1)}{q-q^{-1}}$, whence the first family of symmetries in the statement of Lemma (1). Our above argument thus demonstrates that in the case we consider now ($\pi(\mathsf{k})(z)=\gamma z$, $\pi(\mathsf{e})(z)\ne 0$) there could be no other symmetries. On the other hand, a routine verification shows that those formulas pass through all the relations \eqref{qpr}, \eqref{Lpr}, \eqref{kk1} -- \eqref{effe}, hence determine a family of well-defined $U_q\left(\mathfrak{sl}_2\right)$-symmetries on $\mathbb{C}[z^{\pm 1}]$.

Now consider separately the case $\pi(\mathsf{e})(z)=0$. This means $\pi(\mathsf{e})$ is identically zero on $\mathbb{C}[z^{\pm 1}]$. Apply again \eqref{effe} to $z$ to deduce that now $\gamma=\pm 1$. Of course, $\pi(\mathsf{f})\equiv 0$ in this case, because (a relation similar to) \eqref{rd} fails. We thus obtain two additional symmetries as in the claim of Lemma (1); these could be readily verified to be well defined.

One can verify that the isomorphisms of actions with fixed $\gamma$ and different $a$ are given by $\Psi(z)=\alpha z$ for suitable $\alpha$. The isomorphism of actions corresponding to $\gamma$, $\gamma^{-1}$ is determined by $\Psi(z)=z^{-1}$. This exhausts the action of the group of automorphisms on the space of parameters of symmetries, hence no other isomorphisms between symmetries.

2). Let $\pi(\mathsf{k})(z)=\gamma z^{-1}$.

It is a matter of direct computation that $\pi(\mathsf{k})^2=\operatorname{id}$, whence
$$
\pi(\mathsf{e})=\pi(\mathsf{k^2e})=q^4\pi(\mathsf{ek^2})=q^4\pi(\mathsf{e}),
$$
that is $\pi(\mathsf{e})\equiv 0$. A similar argument proves that $\pi(\mathsf{f})\equiv 0$. Thus \eqref{effe} is satisfied.

The isomorphism between the symmetries with different $\gamma$ is given by $\Psi(z)=\alpha z$ for suitable $\alpha$. \hfill$\blacksquare$

\medskip

\begin{proposition}\label{JB1}
There exist no $U_q\left(\mathfrak{sl}_2\right)$-symmetries on $\mathbb{C}_q[x^{\pm 1},y^{\pm 1}]$ with $\mathsf{k}$ acting via an automorphism $\varphi_{\sigma,\alpha,\beta}$ such that the matrix $\sigma$ has eigenvalues $\lambda=\mu=1$ and a one dimensional eigenspace.
\end{proposition}

{\bf Proof.} Suppose that such symmetry $\pi$ exists. Clearly an eigenvector $\binom{v_1}{v_2}$ of $\sigma$ can be chosen so that $v_1$, $v_2$ are coprime integers. Let $u_1$, $u_2$ be such integers that $u_1v_1+u_2v_2=1$. Consider the matrix $\theta=\begin{pmatrix}v_1 & -u_2\\ v_2 & u_1\end{pmatrix}\in SL(2,\mathbb{Z})$ together with the automorphism $\Phi=\varphi_{\theta,1,1}$ of $\mathbb{C}_q[x^{\pm 1},y^{\pm 1}]$ as in \eqref{Auto}. The isomorphic symmetry $\Phi^{-1}\circ\pi\circ\Phi$ has the $\mathsf{k}$-action as in \eqref{Auto} with the matrix $\theta^{-1}\sigma\theta$, which is of the form $\begin{pmatrix}1 & l\\ 0 & 1\end{pmatrix}$ for some $l\in\mathbb{Z}\setminus\{0\}$. Thus we may assume that $\sigma$ itself has this form.

Now by \eqref{Auto} we have that $\pi(\mathsf{k})(x)=\alpha x$, $\pi(\mathsf{k})(y)=\beta x^ly$ for some $\alpha,\beta\in\mathbb{C}\setminus\{0\}$. Let $\pi(\mathsf{e})(x)=\sum_{i,j}a_{ij}x^iy^j$. A direct computation which uses the relation
$$(x^ry^s)^i=q^{\frac{i(i-1)}2rs}x^{ri}y^{si},\qquad i,r,s\in\mathbb{Z},$$
shows that
\begin{equation}\label{ek(x)}
q^2\pi(\mathsf{ek})(x)=q^2\alpha\sum_{i,j}a_{ij}x^iy^j,
\end{equation}
\begin{equation}\label{ke(x)}
\pi(\mathsf{ke})(x)=
\sum_{i,j}a_{ij}\alpha^i\beta^jq^{\frac{j(j-1)}2\,l}x^{i+lj}y^j.
\end{equation}
Since $l\ne 0$, a comparison of \eqref{ek(x)} and \eqref{ke(x)} demonstrates that if $a_{ij}\ne 0$ for some $i,j$ with $j\ne 0$ then \eqref{ke} fails. So we deduce that $\pi(\mathsf{e})(x)\in\mathbb{C}[x^{\pm 1}]$. Similarly one proves that $\pi(\mathsf{f})(x)\in\mathbb{C}[x^{\pm 1}]$. It follows that the action $\pi$ of $U_q\left(\mathfrak{sl}_2\right)$ leaves invariant the subalgebra $\mathbb{C}[x^{\pm 1}]$, thus defining a $U_q\left(\mathfrak{sl}_2\right)$-symmetry on $\mathbb{C}[x^{\pm 1}]$. By Lemma \ref{Czsym}, we need to consider the two cases.

(A). Let $\alpha^{r-1}=q^2$ for some $r\in\mathbb{Z}$ and $a\in\mathbb{C}\setminus\{0\}$ be such that
$$
\pi(\mathsf{e})(x)=\frac{a}{q^2-1}x^r;\qquad
\pi(\mathsf{f})(x)=q^3(\alpha-1)a^{-1}x^{2-r}.
$$
With $\pi(\mathsf{f})(y)=\sum_{i,j}d_{ij}x^iy^j$ we compute using \eqref{f(z^p)}:
\begin{multline}\label{fk(y)}
q^{-2}\pi(\mathsf{fk})(y)=\beta q^{-2}\pi(\mathsf{f})(x^ly)=\beta q^{-2}\pi(\mathsf{f})(x^l)y+\beta q^{-2}\pi(\mathsf{k})^{-1}(x^l)\pi(\mathsf{f})(y)=
\\ =\beta(\alpha^{-l}-1)qa^{-1}x^{l-r+1}y+
\beta\alpha^{-l}q^{-2}x^l\pi(\mathsf{f})(y)=
\\ =\beta(\alpha^{-l}-1)qa^{-1}x^{l-r+1}y+
\beta\alpha^{-l}q^{-2}\sum_{i,j}d_{ij}x^{i+l}y^j;
\end{multline}
\begin{equation}\label{kf(y)}
\pi(\mathsf{kf})(y)=
\sum_{i,j}d_{ij}\alpha^i\beta^jq^{\frac{j(j-1)}2\,l}x^{i+lj}y^j.
\end{equation}
Since $l\ne 0$, a comparison of \eqref{fk(y)} and \eqref{kf(y)} demonstrates that if $d_{ij}\ne 0$ for some $i,j$ with $j\ne 1$ then \eqref{ke} fails. This implies that
\begin{equation}\label{f(y)x}
\pi(\mathsf{f})(y)x=qx\pi(\mathsf{f})(y).
\end{equation}

Let us use \eqref{f(y)x} and \eqref{f(z^p)} to compute
\begin{multline*}
0=\pi(\mathsf{kf}-q^{-2}\mathsf{fk})(y)=
\pi(\mathsf{kf})(y)-q^{-2}\beta\pi(\mathsf{f})(x^l)y-
q^{-2}\beta\pi(\mathsf{k})^{-1}(x^l)\pi(\mathsf{f})(y)=
\\ =\pi(\mathsf{kf})(y)-
q^{-2}\beta\frac{\alpha^{-l}-1}{\alpha-1}x^{l-1}\pi(\mathsf{f})(x)y- q^{-2}\beta\alpha^{-l}x^l\pi(\mathsf{f})(y),
\end{multline*}
whence
\begin{equation}\label{1s}
\pi(\mathsf{kf})(y)-q^{-2}\alpha^{-l}\beta x^l\pi(\mathsf{f})(y)-
q^{-2}\beta\frac{\alpha^{-l}-1}{\alpha-1}x^{l-1}\pi(\mathsf{f})(x)y=0.
\end{equation}

Furthermore, an application of \eqref{f(y)x} and the explicit form of $\pi(\mathsf{f})(x)$ in the case we consider now yields
\begin{multline*}
0=\pi(\mathsf{f})(yx-qxy)=
\\ =\pi(\mathsf{f})(y)x+\pi(\mathsf{k})^{-1}(y)\pi(\mathsf{f})(x)-
q\pi(\mathsf{f})(x)y-q\pi(\mathsf{k})^{-1}(x)\pi(\mathsf{f})(y)=
\\ =q(1-\alpha^{-1})x\pi(\mathsf{f})(y)+
q^{2-r}\alpha^l\beta^{-1}x^{-l}\pi(\mathsf{f})(x)y-q\pi(\mathsf{f})(x)y,
\end{multline*}
whence
\begin{equation}\label{2'}
\pi(\mathsf{f})(y)+
\frac{\alpha}{\alpha-1}\left(q^{1-r}\alpha^l\beta^{-1}x^{-l-1}-x^{-1}\right)
\pi(\mathsf{f})(x)y=0.
\end{equation}

We need two derived relations. The first one is just $-\pi(\mathsf{k})$ applied to \eqref{2'}:
\begin{equation}\label{2s}
-\pi(\mathsf{kf})(y)+\frac{\alpha q^{-2}}{\alpha-1}\left(\beta x^{l-1}-q^{1-r}x^{-1}\right)\pi(\mathsf{f})(x)y=0.
\end{equation}
The next derived relation is nothing more than \eqref{2'} multiplied on the left by $q^{-2}\alpha^{-l}\beta x^l$:
\begin{equation}\label{3s}
q^{-2}\alpha^{-l}\beta x^l\pi(\mathsf{f})(y)+\frac{\alpha q^{-2}}{\alpha-1}\left(q^{1-r}x^{-1}-\alpha^{-l}\beta x^{l-1}\right)\pi(\mathsf{f})(x)y=0.
\end{equation}

Finally, sum up \eqref{1s}, \eqref{2s}, and \eqref{3s} to obtain
$$
\frac{\beta q^{-2}}{\alpha-1}(1+\alpha)\left(1-\alpha^{-l}\right)x^{l-1}
\pi(\mathsf{f})(x)y=0.
$$
Since $\mathbb{C}_q[x^{\pm 1},y^{\pm 1}]$ is a domain, we conclude that some constant multiplier in this product should be zero. However, in the special case (A) we consider now this can not happen. Thus we obtain a contradiction.

(B). Let $\alpha=\pm 1$, $\pi(\mathsf{e})(x)=\pi(\mathsf{f})(x)=0$.

We have
$$
0=\pi(\mathsf{e})(yx-qxy)=
\pi(\mathsf{e})(y)\pi(\mathsf{k})(x)-qx\pi(\mathsf{e})(y)=
\alpha\pi(\mathsf{e})(y)x-qx\pi(\mathsf{e})(y),
$$
whence
$$\pi(\mathsf{e})(y)x=\alpha^{-1}qx\pi(\mathsf{e})(y).$$
This quasi-commutation relation is possible only if $\pi(\mathsf{e})(y)=\varphi y^p$ for some $\varphi\in\mathbb{C}[x^{\pm 1}]$, $p\in\mathbb{Z}$.

In the case $\alpha=-1$ this can not happen because one should have $q^{p-1}=-1$. The latter implies that $p-1\ne 0$ and $q^{2(p-1)}=1$, which contradicts to our assumptions on $q$.

It remains to assume that $\alpha=1$. In this case $p=1$, and one has that $\pi(\mathsf{e})(y)=\varphi y$, $\pi(\mathsf{k})(\varphi)=\varphi$, hence
$$
0=\pi(\mathsf{ke}-q^2\mathsf{ek})(y)=\pi(\mathsf{k})(\varphi)\beta x^ly-q^2\beta\pi(\mathsf{e})(x^ly)=\beta x^l\varphi y-q^2\beta x^l\varphi y,
$$
whence
$$\beta(1-q^2)x^l\varphi y=0.$$
This implies $\varphi=0$, hence $\pi(\mathsf{e})(y)=0$. Therefore $\pi(\mathsf{e})$ is identically zero on $\mathbb{C}_q[x^{\pm 1},y^{\pm 1}]$. Since $\pi(\mathsf{k})(y)\ne\pi(\mathsf{k})^{-1}(y)$, we observe that \eqref{effe} being applied to $y$ fails. Thus we obtain the final contradiction, which completes the proof of Proposition. \hfill $\blacksquare$

\medskip

(i-2) Suppose that $\operatorname{tr}\sigma=1$. This means that
$\lambda=\frac12+i\frac{\sqrt{3}}2$, $\mu=\frac12-i\frac{\sqrt{3}}2$. In
particular, the matrix $\sigma$ has a finite order, more precisely
$\sigma^6=I$. Hence the same is true for the corresponding automorphism of
$\mathbb{C}_q[x^{\pm 1},y^{\pm 1}]$ as in \eqref{Auto} with
$\alpha=\beta=1$. However, we need a more subtle claim.

\begin{lemma}\label{fo}
Assume we are given an arbitrary pair $(\alpha,\beta)\in(\mathbb{C}^*)^2$
and a matrix $\sigma\in SL(2,\mathbb{Z})$ with the properties listed in the
subcase (ii-2) as well as also in the subcases (ii-3), (ii-4) below. Then
the automorphism $\varphi_{\sigma,\alpha,\beta}$ of $\mathbb{C}_q[x^{\pm
1},y^{\pm 1}]$ determined by \eqref{Auto} has a finite order, the latter
being larger than $2$.
\end{lemma}

{\bf Proof.} One readily computes that $\det(\sigma-I)=2-\operatorname{tr}\sigma=1$, hence the inverse matrix $(\sigma-I)^{-1}$ is integral. Thus we have a well-defined pair $(\alpha',\beta')=(\sigma-I)(\alpha,\beta)(\sigma-I)^{-1}\in(\mathbb{C}^*)^2$ as in \eqref{sp}. Now a simple computation shows that in the group $\operatorname{Aut}(\mathbb{C}_q[x^{\pm 1},y^{\pm 1}])$ one has the conjugation $$(\alpha',\beta')^{-1}\sigma(\alpha',\beta')=(\alpha,\beta)\sigma.$$ It follows that $\varphi_{\sigma,\alpha,\beta}^6=\operatorname{id}$, which proves Lemma in the present subcase (i-2). It will become clear below that in the subcases (i-3), (i-4) this proof requires only minor modifications. \hfill $\blacksquare$

\medskip

Now assume that $\pi$ is a symmetry with $\sigma$ possessing the properties
listed in the subcase (i-2). It follows from \eqref{ke} that
$\mathsf{k}^6\mathsf{e}\mathsf{k}^{-6}=q^{12}\mathsf{e}$. By Lemma
\ref{fo}, $\pi(\mathsf{k})^6=\operatorname{id}$, hence
$q^{12}\pi(\mathsf{e})=\pi(\mathsf{e})$. Since $q$ is not a root of $1$,
this implies that $\pi(\mathsf{e})$ is identically zero. A similar argument
establishes also that $\pi(\mathsf{f})\equiv 0$.

Now let us apply $\pi$ to \eqref{effe}. The above observations show that we
obtain identically zero in the l.h.s. But this is not the case with the
r.h.s, because $\pi(\mathsf{k})\not\equiv\pi(\mathsf{k})^{-1}$, as by Lemma
\ref{fo} the order of $\pi(\mathsf{k})$ is larger than $2$. The
contradiction we obtain this way proves that {\it there exist no symmetries
corresponding to $\sigma$ as in the subcase (i-2)}.

\medskip

(i-3) Suppose that $\operatorname{tr}\sigma=0$. This means that
$\lambda=i$, $\mu=-i$. The matrix $\sigma$ has order $4$, $\sigma^4=I$,
hence $\varphi_{\sigma,1,1}^4=\operatorname{id}$.

To prove Lemma \ref{fo} in this subcase, observe first that the subgroup
$T\subset(\mathbb{C}^*)^2\subset\operatorname{Aut}(\mathbb{C}_q[x^{\pm
1},y^{\pm 1}])$ formed by $(r,s)$ with $r,s=\pm 1$, is normal in
$\operatorname{Aut}(\mathbb{C}_q[x^{\pm 1},y^{\pm 1}])$, as one can see
from \eqref{sp}. Therefore, the subgroup of
$\operatorname{Aut}(\mathbb{C}_q[x^{\pm 1},y^{\pm 1}])$ generated by $T$
and $\sigma^i$, $i=0,1,2,3$, is finite.

In this subcase one has $\det(\sigma-I)=2-\operatorname{tr}\sigma=2$, hence
$(\sigma-I)\sigma'=2I$ for some integral matrix $\sigma'$. As (the second
equality of) \eqref{sp} determines a well-defined action of the semigroup
of integral matrices on $(\mathbb{C}^*)^2$, we conclude that, given
arbitrary $(\alpha,\beta)\in(\mathbb{C}^*)^2$, one has
$$
(\alpha,\beta)^{(\sigma-I)\sigma'}=
\left((\alpha,\beta)^{\sigma'}\right)^{(\sigma-I)}=(\alpha,\beta)^2.
$$
After passage to square roots we obtain
$$
\left((\alpha',\beta')^{\sigma'}\right)^{(\sigma-I)}=(r,s)(\alpha,\beta),
$$
for some $\alpha'$, $\beta'$, $r$, $s$ such that $\alpha'^2=\alpha$,
$\beta'^2=\beta$; $r,s\in\{-1,1\}$. Now set
$(\alpha'',\beta'')\stackrel{\operatorname{def}}{=}
(\alpha',\beta')^{\sigma'}$. A routine computation shows that one has the
conjugation
$$
(\alpha'',\beta'')^{-1}(r,s)\sigma(\alpha'',\beta'')=(\alpha,\beta)\sigma.
$$
By the above argument, $(r,s)\sigma$ has finite order as an element of a finite subgroup, hence the same is true for its conjugate $(\alpha,\beta)\sigma$. This order can not be less than the order of the projection $\sigma$ to $SL(2,\mathbb{Z})$, which is $4$. This proves Lemma \ref{fo}. It remains to proceed the same way as in the subcase (i-2) to conclude that {\it there exist no symmetries corresponding to $\sigma$ as in the subcase (i-3)}.

\medskip

(i-4) Suppose that $\operatorname{tr}\sigma=-1$. This means that
$\lambda=-\frac12+i\frac{\sqrt{3}}2$, $\mu=-\frac12-i\frac{\sqrt{3}}2$. The
matrix $\sigma$ has order $3$, $\sigma^3=I$, hence
$\varphi_{\sigma,1,1}^3=\operatorname{id}$.

To prove Lemma \ref{fo} in this subcase, observe first that the subgroup
$T\subset(\mathbb{C}^*)^2\subset\operatorname{Aut}(\mathbb{C}_q[x^{\pm
1},y^{\pm 1}])$ formed by $(r,s)$ with $r,s=\zeta^i$,
$\zeta=-\frac12+i\frac{\sqrt{3}}2$, $i=0,1,2$, is normal in
$\operatorname{Aut}(\mathbb{C}_q[x^{\pm 1},y^{\pm 1}])$, as one can see
from \eqref{sp}. Therefore, the subgroup of
$\operatorname{Aut}(\mathbb{C}_q[x^{\pm 1},y^{\pm 1}])$ generated by $T$
and $\sigma^i$, $i=0,1,2$, is finite.

In this subcase one has $\det(\sigma-I)=2-\operatorname{tr}\sigma=3$, hence
$(\sigma-I)\sigma'=3I$ for some integral matrix $\sigma'$. As (the second
equality of) \eqref{sp} determines a well-defined action of the semigroup
of integral matrices on $(\mathbb{C}^*)^2$, we conclude that, given
arbitrary $(\alpha,\beta)\in(\mathbb{C}^*)^2$, one has
$$
(\alpha,\beta)^{(\sigma-I)\sigma'}=
\left((\alpha,\beta)^{\sigma'}\right)^{(\sigma-I)}=(\alpha,\beta)^3.
$$
After passage to cubic roots we obtain
$$
\left((\alpha',\beta')^{\sigma'}\right)^{(\sigma-I)}=(r,s)(\alpha,\beta),
$$
for some $\alpha'$, $\beta'$, $r$, $s$ such that $\alpha'^3=\alpha$,
$\beta'^3=\beta$; $r,s\in\{\zeta^i,\:i=0,1,2\}$. Now set
$(\alpha'',\beta'')\stackrel{\operatorname{def}}{=}
(\alpha',\beta')^{\sigma'}$. A routine computation shows that one has the
conjugation
$$
(\alpha'',\beta'')^{-1}(r,s)\sigma(\alpha'',\beta'')=(\alpha,\beta)\sigma.
$$
By the above argument, $(r,s)\sigma$ has finite order as an element of a finite subgroup, hence the same is true for its conjugate $(\alpha,\beta)\sigma$. This order can not be less than the order of the projection $\sigma$ to $SL(2,\mathbb{Z})$, which is $3$. This proves Lemma \ref{fo}. It remains to proceed the same way as in the subcase (i-2) to conclude that {\it there exist no symmetries corresponding to $\sigma$ as in the subcase (i-4)}.

\medskip

(i-5) Suppose that $\operatorname{tr}\sigma=-2$. This means that
$\lambda=\mu=-1$, hence either $\sigma=-I$ or $\sigma$ is a conjugate matrix to a Jordan block, that is, the eigenspace of $\sigma$ is one dimensional.

\begin{theorem}\label{sigma=-I}
There exists a two-parameter ($\alpha,\beta\in\mathbb{C}^*$) family of
$U_q(\mathfrak{sl}_2)$-symmetries on $\mathbb{C}_q[x^{\pm 1},y^{\pm 1}]$
that correspond to $\sigma=-I$
\begin{align}
\pi(\mathsf{k})(x) &=\alpha^{-1}x^{-1}; & \pi(\mathsf{k})(y) &=\beta^{-1}y^{-1};
& \label{k-I}
\\ \pi(\mathsf{e})(x) &=0;  & \pi(\mathsf{e})(y) &=0; & \label{e-I}
\\ \pi(\mathsf{f})(x) &=0;  & \pi(\mathsf{f})(y) &=0. & \label{f-I}
\end{align}
These are all the symmetries with $\sigma=-I$. These symmetries are all
isomorphic, in particular to that with $\alpha=\beta=1$.
\end{theorem}

{\bf Proof.} A routine verification demonstrates that the action given by
\eqref{k-I} -- \eqref{f-I} passes through all the relations \eqref{qpr},
\eqref{Lpr} in $\mathbb{C}_q[x^{\pm 1},y^{\pm 1}]$ and \eqref{kk1} --
\eqref{eps} in $U_q(\mathfrak{sl}_2)$. This means that this action is
really a $U_q(\mathfrak{sl}_2)$-symmetry.

To see that there are no more symmetries with $\sigma=-I$, observe that
$\alpha$ and $\beta$ are arbitrary non-zero complex numbers, hence
\eqref{k-I} exhausts all the possibilities for the action of $\mathsf{k}$
by an automorphism (see \eqref{Auto}). As for the action of $\mathsf{e}$
and $\mathsf{f}$, note first that $\pi(\mathsf{k})^2=\operatorname{id}$ (a
straightforward computation). Hence by \eqref{ke} one has
$\pi(\mathsf{e})=\pi(\mathsf{k}^2\mathsf{e}\mathsf{k}^{-2})=q^4\pi(\mathsf{e})$.
Since $q$ is not a root of $1$, we deduce that $\pi(\mathsf{e})\equiv 0$.
Similarly, $\pi(\mathsf{f})\equiv 0$.

Let us verify that all the symmetries listed in the Theorem are isomorphic. Given any such symmetry with $\pi(\mathsf{k})$ being $(\alpha,\beta)(-I)\in\operatorname{Aut}(\mathbb{C}_q[x^{\pm 1},y^{\pm 1}])$ and an arbitrary automorphism $A=(\mu,\nu)\tau$, $\tau\in SL(2,\mathbb{Z})$, $(\mu,\nu)\in(\mathbb{C}^*)^2$, one readily computes that
$$A(\alpha,\beta)(-I)A^{-1}=(\mu^2,\nu^2)(\alpha,\beta)^\tau(-I).$$ Clearly, an appropriate choice of $(\mu,\nu)$, $\tau$ makes this conjugate automorphism $(-I)$. \hfill$\blacksquare$

\medskip

\begin{remark}
Although the action of $\mathsf{k}$ in the symmetries of Theorem \ref{sigma=-I} does not reduce to multiplying the generators $x$, $y$ by weight constants as in \cite{DS}, the associated $U_q(\mathfrak{sl}_2)$-actions are weight modules. Namely, a basis of weight vectors in $\mathbb{C}_q[x^{\pm 1},y^{\pm 1}]$ is given by
$$
\{1\}\cup\{u_{ij}=\alpha^i\beta^jx^iy^j+x^{-i}y^{-j}|\:i,j>0\}\cup
\{v_{ij}=\alpha^i\beta^jx^iy^j-x^{-i}y^{-j}|\:i,j>0\},
$$
so that $\pi(\mathsf{k})(1)=1$, $\pi(\mathsf{k})(u_{ij})=u_{ij}$, $\pi(\mathsf{k})(v_{ij})=-v_{ij}$.
\end{remark}

\medskip

\begin{proposition}
There exist no $U_q\left(\mathfrak{sl}_2\right)$-symmetries on $\mathbb{C}_q[x^{\pm 1},y^{\pm 1}]$ with $\mathsf{k}$ acting via an automorphism $\varphi_{\sigma,\alpha,\beta}$ such that the matrix $\sigma$ has eigenvalues $\lambda=\mu=-1$ and a one dimensional eigenspace.
\end{proposition}

P r o o f . Suppose that such symmetry $\pi$ exists. One can readily apply the same argument as that at the beginning of the proof of Proposition \ref{JB1} in order to reduce matters to the case $\sigma=\begin{pmatrix}-1 & l\\ 0 & -1\end{pmatrix}$ for some $l\in\mathbb{Z}\setminus\{0\}$.

Now by \eqref{Auto} we have that $\pi(\mathsf{k})(x)=\alpha x^{-1}$, $\pi(\mathsf{k})(y)=\beta x^ly^{-1}$ for some $\alpha,\beta\in\mathbb{C}\setminus\{0\}$. Let $\pi(\mathsf{e})(x)=\sum_{i,j}a_{ij}x^iy^j$. A direct computation which uses \eqref{ext_e} shows that
\begin{equation}\label{ek(x)_}
q^2\pi(\mathsf{ek})(x)=-\sum_{i,j}q^{j+2}a_{ij}x^iy^j,
\end{equation}
\begin{equation}\label{ke(x)_}
\pi(\mathsf{ke})(x)=
\sum_{i,j}a_{ij}\alpha^i\beta^jq^{-\frac{j(j-1)}2\,l}x^{-i+lj}y^{-j}.
\end{equation}
Let us compare \eqref{ek(x)_}, \eqref{ke(x)_}. Equate the coefficients at $x^i$, and $x^{-i}$, respectively, $i\in\mathbb{Z}$, $j=0$:
$$-q^2a_{i0}=a_{-i,0}\alpha^{-i},\qquad -q^2a_{-i,0}=a_{i0}\alpha^i.$$
The product of two latter relations is $q^4a_{i,0}a_{-i,0}=a_{i,0}a_{-i,0}$. This clearly implies that $a_{i,0}a_{-i,0}=0$. Due to the above relations, it follows that $a_{i,0}=0$ for all $i$.

Now let $j\ne 0$. Let us compare again \eqref{ek(x)_}, \eqref{ke(x)_} and equate the coefficients at $x^{-i+lj}y^{-j}$ to obtain
$$a_{-i+lj,\,-j}=-q^{-\frac{j(j-1)}2\,l+j-2}\alpha^i\beta^ja_{ij}.$$
One more iteration of this relation yields
$$
a_{i-2lj,\,j}=
-q^{-\frac{(-j)(-j-1)}2\,l-j-2}\alpha^{-i+lj}\beta^{-j}a_{-i+lj,\,-j}=
q^{-j^2l-4}\alpha^{lj}a_{ij}.
$$
Since $l\ne 0$, this implies that, once some $a_{ij}$ with $j\ne 0$ is assumed to be non-zero, then infinitely many of other $a_{ij}$'s appear to be non-zero. This is certainly impossible, because the element $\pi(\mathsf{e})(x)=\sum_{i,j}a_{ij}x^iy^j$ of the algebra $\mathbb{C}_q[x^{\pm 1},y^{\pm 1}]$ is anyway a finite sum of monomials. This proves that with $j\ne 0$, $a_{ij}=0$ for all $i\in\mathbb{Z}$.

Thus we conclude that $\pi(\mathsf{e})(x)=0$. Now let us compute
$$
0=\pi(\mathsf{e})(yx-qxy)=
\alpha\pi(\mathsf{e})(y)x^{-1}-qx\pi(\mathsf{e})(y),
$$
whence $\alpha\pi(\mathsf{e})(y)x^{-1}=qx\pi(\mathsf{e})(y)$. The latter can not happen with non-zero $\pi(\mathsf{e})(y)$ because the maximum degree in $x$ of monomials in the l.h.s. is lower than such degree in the r.h.s. Therefore $\pi(\mathsf{e})(y)=0$, hence $\pi(\mathsf{e})$ is identically zero on $\mathbb{C}_q[x^{\pm 1},y^{\pm 1}]$. On the other hand, $\pi(\mathsf{k})(y)\ne\pi(\mathsf{k})^{-1}(y)$, hence \eqref{effe} being applied to $y$ fails. This contradiction completes the proof of Proposition. \hfill $\blacksquare$

\medskip

Let us consider the case (ii).

\begin{proposition}\label{real_eigens}
There exist no $U_q(\mathfrak{sl}_2)$-symmetries on $\mathbb{C}_q[x^{\pm
1},y^{\pm 1}]$ with $\mathsf{k}$ acting via an automorphism
$\varphi_{\sigma,\alpha,\beta}$ such that the matrix $\sigma$ has
eigenvalues $\lambda,\mu\in\mathbb{R}\setminus\{-1,1\}$.
\end{proposition}

To prove this Proposition, we need some observations. Firstly, given a symmetry $\pi$, let us introduce the notation
\begin{equation}\label{exey}
\pi(\mathsf{e})(x)=\sum_{i,j}a_{ij}x^iy^j,\qquad
\pi(\mathsf{e})(y)=\sum_{i,j}b_{ij}x^iy^j.
\end{equation}
Here, of course, only finitely many of $a_{ij}$, $b_{ij}$ are non-zero. Let
us denote $D=\left\{\binom{i}{j}\in\mathbb{Z}^2:\:a_{ij}\ne 0\right\}$,
$E=\left\{\binom{i}{j}\in\mathbb{Z}^2:\:b_{ij}\ne 0\right\}$.

A straightforward induction argument allows one to deduce the relation
\begin{equation}\label{x^p>0}
\pi(\mathsf{e})(x^p)=
\sum_{r=0}^{p-1}x^{p-1-r}\pi(\mathsf{e})(x)\pi(\mathsf{k})(x)^r,
\end{equation}
for integers $p>0$, together with a similar but slightly different formula
for $p<0$ (just due to \eqref{ext_e}):
\begin{equation}\label{x^p<0}
\pi(\mathsf{e})(x^p)=
-\sum_{r=0}^{-p-1}x^{p+r}\pi(\mathsf{e})(x)\pi(\mathsf{k})(x)^{-r-1},
\end{equation}

One can combine these two formulas in order to obtain a universal formula
for all $p\in\mathbb{Z}\setminus\{0\}$, which, however, ignores the specific values of constant multipliers at monomials. More precisely, we have
\begin{equation}\label{e(x^p)_univ}
\pi(\mathsf{e})(x^p)=
\sum_{\binom{i}{j}\in D}\sum_{r=\min\{0,p\}}^{\max\{0,p\}-1}
\operatorname{const}(i,j,p,r,k,m)x^{p-1-r+i+rk}y^{j+rm},
\end{equation}
where $k$ and $m$ determine the action of $\mathsf{k}$ as in \eqref{Auto},
the constants $a_{ij}$ are as in \eqref{exey}, and all the constants
$\operatorname{const}(\ldots)$ being non-zero. This formula is going to be
useful in the sequel where we are about to compute the sums of monomials
modulo constant multipliers.

In a similar way, one has
\begin{equation}\label{y^p>0}
\pi(\mathsf{e})(y^p)=
\sum_{r=0}^{p-1}y^{p-1-r}(\pi(\mathsf{e})y)(\pi(\mathsf{k})y)^r,
\end{equation}
for integers $p>0$, and
\begin{equation}\label{y^p<0}
\pi(\mathsf{e})(y^p)=
-\sum_{r=0}^{-p-1}y^{p+r}(\pi(\mathsf{e})y)(\pi(\mathsf{k})y)^{-r-1},
\end{equation}
for $p<0$.

Again, the latter two formulas can be combined in order to obtain a
universal formula for all $p\in\mathbb{Z}\setminus\{0\}$, similar to
\eqref{e(x^p)_univ}:
\begin{equation}\label{e(y^p)_univ}
\pi(\mathsf{e})(y^p)=
\sum_{\binom{i}{j}\in E}\sum_{r=\min\{0,p\}}^{\max\{0,p\}-1}
\operatorname{const}(i,j,p,r,l,n)x^{i+rl}y^{p-1-r+j+rn},
\end{equation}
where $l$ and $n$ determine the action of $\mathsf{k}$ as in \eqref{Auto},
the constants $b_{ij}$ are as in \eqref{exey}, and all the constants
$\operatorname{const}(\ldots)$ being non-zero.

Another observation is related to the specific form of the (integral)
entries of the matrix
$$\sigma^N=\begin{pmatrix}a(N) & b(N)\\ c(N) & d(N)\end{pmatrix}.$$

Let $\lambda$, $\lambda^{-1}$ be the real eigenvalues of $\sigma$ as under
the assumptions of Proposition \ref{real_eigens}, $|\lambda|>1$, then
$$
\sigma^N=
\Phi\begin{pmatrix}\lambda^N & 0\\ 0 & \lambda^{-N}\end{pmatrix}\Phi^{-1}
$$
for some invertible complex matrix $\Phi$ and $N\in\mathbb{Z}$. A routine
computation shows that
$$
\sigma^N=\begin{pmatrix}a\lambda^N+a'\lambda^{-N} & b\lambda^N+b'\lambda^{-N}
\\ c\lambda^N+c'\lambda^{-N} & d\lambda^N+d'\lambda^{-N}\end{pmatrix},
$$
with some $a,a',b,b',c,c',d,d'\in\mathbb{C}$. Substitute now $N=0$; as
$\sigma^0=I$, we deduce that in fact
$$
\sigma^N=\begin{pmatrix}a\lambda^N+(1-a)\lambda^{-N} &
b(\lambda^N-\lambda^{-N})\\ c(\lambda^N-\lambda^{-N}) &
d\lambda^N+(1-d)\lambda^{-N}\end{pmatrix}.
$$
Finally, computing $\det\sigma^N$, which is just $1$, we find that
$ad-bc=0$ and $d=1-a$.

Note that under the assumptions of Proposition \ref{real_eigens} $b\ne 0$. In fact, if $b=0$, then $ad=0$, hence either $a=0$ or $a=1$. In both cases $\sigma$ becomes a triangular matrix whose diagonal entries are $\lambda$, $\lambda^{-1}$. Since these are integers, one has $\lambda=\pm 1$, which is not our case.

In a similar way one deduces that $c\ne 0$, $a\ne 0$, $a\ne 1$.

Thus we have $c=\frac{a(1-a)}{b}$, so
\begin{equation}\label{sigma^N}
\sigma^N=\begin{pmatrix}a(N) & b(N)\\ c(N) & d(N)\end{pmatrix}=
\begin{pmatrix}a\lambda^N+(1-a)\lambda^{-N} &
b(\lambda^N-\lambda^{-N})\\ \frac{a(1-a)}{b}(\lambda^N-\lambda^{-N}) &
(1-a)\lambda^N+a\lambda^{-N}\end{pmatrix}.
\end{equation}

\medskip

{\bf Proof} of Proposition \ref{real_eigens}. Assume the contrary, that is, there exists a symmetry $\pi$ with the properties as in Proposition \ref{real_eigens}. We restrict our considerations to the case when $\mathsf{k}$ acts via an automorphism $\varphi_{\sigma,\alpha,\beta}$ with $\alpha=\beta=1$. It will become clear in what follows that this extra assumption is not really restrictive. However, this restriction will be implicit unless the contrary is stated explicitly.

Let us use \eqref{Auto}, \eqref{def}, \eqref{e(x^p)_univ},
\eqref{e(y^p)_univ}, to compute (modulo non-zero constant multipliers at monomials)
\begin{multline}\label{ekN(x)}
\pi\left(\mathsf{e}\mathsf{k}^N\right)(x)=
\pi(\mathsf{e})\left(x^{a(N)}y^{c(N)}\right)=
\\ =x^{a(N)}\pi(\mathsf{e})\left(y^{c(N)}\right)+
\pi(\mathsf{e})\left(x^{a(N)}\right)\pi(\mathsf{k})\left(y^{c(N)}\right)=
\\ =x^{a(N)}\sum_{\binom{i}{j}\in E}\sum_{r=\min\{0,c(N)\}}^{\max\{0,c(N)\}-1}
\operatorname{const} x^{i+rb(1)}y^{j+c(N)-1+r(d(1)-1)}+
\\ +\sum_{\binom{i}{j}\in D}\sum_{r=\min\{0,a(N)\}}^{\max\{0,a(N)\}-1}
\operatorname{const}x^{i+a(N)-1+r(a(1)-1)}y^{j+rc(1)}
\left(x^{b(1)}y^{d(1)}\right)^{c(N)}=
\\ =\sum_{\binom{i}{j}\in E}\sum_{r=\min\{0,c(N)\}}^{\max\{0,c(N)\}-1}
\operatorname{const}x^{i+a(N)+rb(1)}y^{j+c(N)-1+r(d(1)-1)}+
\\ +\sum_{\binom{i}{j}\in D}\sum_{r=\min\{0,a(N)\}}^{\max\{0,a(N)\}-1}
\operatorname{const}x^{i+a(N)+b(1)c(N)-1+r(a(1)-1)}y^{j+d(1)c(N)+rc(1)}.
\end{multline}
The constant multipliers here are all non-zero. In the case when either $a(N)=0$ or $c(N)=0$, the corresponding sum in \eqref{ekN(x)} is totally
absent (equals to $0$), which agrees with $\pi(\mathsf{e})(\mathbf{1})=0$. On the other hand,
\begin{multline}\label{kNe(x)}
\pi\left(\mathsf{k}^N\mathsf{e}\right)(x)=\pi\left(\mathsf{k}^N\right)
\left(\sum_{\binom{i}{j}\in D}\operatorname{const}x^iy^j\right)=
\\ =\sum_{\binom{i}{j}\in D}\operatorname{const}
\left(x^{a(N)}y^{c(N)}\right)^i\left(x^{b(N)}y^{d(N)}\right)^j=
\\ =\sum_{\binom{i}{j}\in D}\operatorname{const}
x^{ia(N)+jb(N)}y^{ic(N)+jd(N)},
\end{multline}
with all the constant multipliers being non-zero.

Now let $\xi=\binom{\xi_1}{\xi_2}$, $\eta=\binom{\eta_1}{\eta_2}$ be the
eigenvectors of $\sigma$ corresponding to the eigenvalues $\lambda$ and
$\lambda^{-1}$, with $|\lambda|>1$ and $|\lambda^{-1}|<1$, respectively.
One may assume that $\xi_i$, $\eta_i$ are real, so that everything is
embedded into the real vector space $\mathbb{R}^2$, together with the
action of $SL(2,\mathbb{Z})$ on it, which leaves invariant the integral
lattice $\mathbb{Z}^2$.

Note that $\dfrac{\eta_2}{\eta_1}$ is irrational. In fact, if one assumes the contrary, one can normalize $\eta$ so that it has integral coordinates, together with all $\sigma^N\eta$, $N\in\mathbb{Z}$. Since the latter vectors are just $\lambda^{-N}\eta$, this makes a contradiction. In a similar way one observes that $\dfrac{\xi_2}{\xi_1}$ is irrational.

Consider linear functionals $\Phi$, $\Psi$ on the real vector space
$\mathbb{R}^2$ given by
$$\Phi(s\xi+t\eta)=s,\quad\Psi(s\xi+t\eta)=t,\qquad s,t\in\mathbb{R}.$$
Set also
$$L_\varepsilon=\left\{w\in\mathbb{R}^2:\:|\Psi(w)|<\varepsilon\right\},$$
for $\varepsilon>0$.

We are going to use the one-to-one correspondence
$\operatorname{const}x^iy^j\mapsto\binom{i}{j}$ between the monomials in
$\mathbb{C}_q[x^{\pm 1},y^{\pm 1}]$ modulo non-zero constant multipliers
and pairs of integers in $\mathbb{Z}^2$. In the latter picture, one should
certainly expect, as a consequence of \eqref{ke}, that the (finite)
collections of pairs of integers coming from \eqref{ekN(x)} and
\eqref{kNe(x)} coincide. But the collections themselves look very
different. As for \eqref{kNe(x)}, the corresponding subset in
$\mathbb{Z}^2\subset\mathbb{R}^2$ is nothing more than $\sigma^ND$, which
is inside a narrow stripe $L_\varepsilon$, $\varepsilon>0$ being as small
as desired with $N\in\mathbb{Z}_+$ big enough. On the other hand, the set
of pairs coming from \eqref{ekN(x)} is formed by the two collections of
arithmetic progressions
\begin{multline}\label{pE}
\left\{\binom{i+a(N)+rb(1)}{j+c(N)-1+r(d(1)-1)}:\right.
\\ \left.\phantom{\binom{N}{N}}r\in\mathbb{Z},\quad\min\{0,c(N)\}\le r<\max\{0,c(N)\}\right\},
\end{multline}
with $\binom{i}{j}\in E$, and
\begin{multline}\label{pD}
\left\{\binom{i+a(N)+b(1)c(N)-1+r(a(1)-1)}{j+d(1)c(N)+rc(1)}:\right.
\\ \left.\phantom{\binom{N}{N}}r\in\mathbb{Z},\quad\min\{0,a(N)\}\le r<\max\{0,a(N)\}\right\},
\end{multline}
with $\binom{i}{j}\in D$. In the case when either $a(N)=0$ or $c(N)=0$, the corresponding set of progressions is treated as void because the associated sum in \eqref{ekN(x)} is totally absent (equals to $0$). The steps of these progressions, $h_E=\binom{b(1)}{d(1)-1}$ and $h_D=\binom{a(1)-1}{c(1)}$, respectively, do not depend on $N$. These vectors are certainly linear independent, being the columns of the non-degenerate matrix $\sigma-I$.

So, the only problem in producing a desired contradiction is in observing
that some pairs in \eqref{pE}, \eqref{pD} may fail to be distinct, hence
the corresponding monomials in \eqref{ekN(x)} may fail to survive after
possible reductions.

\begin{lemma}\label{2pt}
Each of the sets $D$ and $E$ contains at most two elements.
\end{lemma}

{\bf Proof.} We prove this Lemma for $D$. A completely similar argument can be used to prove it for $E$.

The above argument on irrationality implies that all the values
$\Phi\binom{i}{j}$, $\binom{i}{j}\in D$, are pairwise different, hence
$$
d_\Phi=\min\left\{\left|\Phi\binom{i}{j}-\Phi\binom{i'}{j'}\right|:\;
\binom{i}{j}\ne\binom{i'}{j'}\in D\right\}>0,
$$
and
$$
|\lambda^N|d_\Phi=
\min\left\{\left|\Phi\binom{i}{j}-\Phi\binom{i'}{j'}\right|:\;
\binom{i}{j}\ne\binom{i'}{j'}\in\sigma^ND\right\}
$$
for $N\in\mathbb{Z}$. Of course, the latter inequality ($d_\Phi>0$)
anticipates $\operatorname{card}D>1$, while the opposite assumption already
implies the claim of Lemma.

In a similar way set
$$
d_\Psi=\max\left\{\left|\Psi\binom{i}{j}-\Psi\binom{i'}{j'}\right|:\;
\binom{i}{j},\:\binom{i'}{j'}\in D\right\}.
$$

Let
$A=\max\left\{\left|\Phi\binom{i}{j}-\Phi\binom{i'}{j'}\right|:\:
\binom{i}{j},\,\binom{i'}{j'}\in D\right\}$. As $\Psi(h_D)\ne 0$, one can
choose $\varepsilon>0$ so that $\varepsilon<\frac12|\Psi(h_D)|$, and then
choose $N>0$ so that $\sigma^ND\subset L_\varepsilon$ and
\begin{equation}\label{lNdphi}
|\lambda^N|d_\Phi>
A+(d_\Psi+2\varepsilon)\left|\frac{\Phi(h_D)}{\Psi(h_D)}\right|+1.
\end{equation}

Consider the second sum in \eqref{ekN(x)} corresponding to $\binom{i}{j}\in
D$, together with the associated set of progressions \eqref{pD} in
$\mathbb{Z}^2$. By our choice of $\varepsilon$, every such progression has
at most one intersection point with $L_\varepsilon$, hence also with
$\sigma^ND$. We claim that in fact all the progressions in \eqref{pD}
together produce at most one intersection point with $\sigma^ND$. To see
this, assume the contrary, that is two progressions as in \eqref{pD} meet
$\sigma^ND$ in two different points. As the steps of these progressions are
the same, these progressions are disjoint. To be more precise, we have some
$\binom{i_1}{j_1},\,\binom{i_2}{j_2}\in\binom{a(N)+b(1)c(N)-1}{d(1)c(N)}+D$
and some $r_1,\,r_2\in\mathbb{Z}$ with $\min\{0,a(N)\}\le
r<\max\{0,a(N)\}$, such that
$$\binom{i_1}{j_1}+r_1h_D,\,\binom{i_2}{j_2}+r_2h_D\in\sigma^ND.$$
In these settings we have
\begin{multline*}
d_\Psi\ge\left|\Psi\binom{i_1}{j_1}-\Psi\binom{i_2}{j_2}\right|\ge
\\ \ge|r_1-r_2||\Psi(h_D)|-\left|\Psi\left(\binom{i_1}{j_1}+r_1h_D\right)-
\Psi\left(\binom{i_2}{j_2}+r_2h_D\right)\right|\ge
\\ \ge|r_1-r_2||\Psi(h_D)|-2\varepsilon,
\end{multline*}
whence
$$|r_1-r_2|\le\frac{d_\Psi+2\varepsilon}{|\Psi(h_D)|}.$$
This estimate, together with \eqref{lNdphi}, implies
\begin{multline*}
A\ge\left|\Phi\binom{i_1}{j_1}-\Phi\binom{i_2}{j_2}\right|\ge
\\ \ge\left|\Phi\left(\binom{i_1}{j_1}+r_1h_D\right)-
\Phi\left(\binom{i_2}{j_2}+r_2h_D\right)\right|-|r_1-r_2||\Phi(h_D)|\ge
\\ \ge|\lambda^N|d_\Phi-\frac{d_\Psi+2\varepsilon}{|\Psi(h_D)|}|\Phi(h_D)|
>A+1.
\end{multline*}
We thus obtain a contradiction which proves the existence of at most one
intersection point of progressions \eqref{pD} and $\sigma^ND$ for $N$
chosen above and all bigger $N$.

A similar argument, possibly after decreasing $\varepsilon$ and increasing
$N$, allows one to establish the existence of additionally at most one
intersection point of progressions \eqref{pE} and $\sigma^ND$. This proves
the claim of Lemma, because, due to \eqref{ke}, the union of points in
\eqref{pD} and \eqref{pE} must contain $\sigma^ND$. \hfill$\blacksquare$

\medskip

\begin{lemma}\label{2tf}
$\operatorname{card}D=\operatorname{card}E=2$, and $D=\left\{\binom{i}{j};\binom{i}{j}+h_D\right\}$, $E=\left\{\binom{i'}{j'};\binom{i'}{j'}+h_E\right\}$ for some integers $i$, $j$, $i'$, $j'$.
\end{lemma}

{\bf Proof.} In view of Lemma \ref{2pt} we have to consider finitely many cases to be discarded.

(a) Let us suppose $\operatorname{card}D=0$. This means that
$\pi(\mathsf{e})(x)=0$. If one also has $\operatorname{card}E=0$, i.e.,
$\pi(\mathsf{e})(y)=0$, then $\pi(\mathsf{e})$ is identically zero. Let us
apply \eqref{effe} via the action $\pi$ to $x$. As the left hand side is
zero, the right hand side being zero is equivalent to
$\pi(\mathsf{k})^2(x)=x$. This already implies (even without assuming
$\alpha=\beta=1$ as at the beginning of the proof) that the matrix
$\sigma^2$ has eigenvalue $1$, which is not our case. Thus we get a
contradiction.

Suppose now in this case $\operatorname{card}E=1$, that is, \eqref{pE} contains just one progression (note that $c(N)\ne 0$), while \eqref{pD} is void. As the corresponding monomials in \eqref{ekN(x)} are linearly independent, we conclude that \eqref{ekN(x)} is non-zero, while \eqref{kNe(x)} is zero. This contradicts to \eqref{ke}.

The last subcase here is $\operatorname{card}E=2$. Again as $c(N)\ne 0$, \eqref{pE} contains just two progressions. These can not coincide as sets of points, as they correspond to different points of $E$ and have the same step $h_E$. There exists a point that belongs to one of the progressions but not to the other. The corresponding monomial in \eqref{ekN(x)} will survive under possible reductions, thus making \eqref{ekN(x)} non-zero, while \eqref{kNe(x)} is zero. This contradicts to \eqref{ke}.

(b) Consider the case $\operatorname{card}D=1$. There exists a single progression in \eqref{pD}. Let us choose $\varepsilon>0$ and $N>0$ as in the proof of Lemma \ref{2pt}, which provides that the progression in \eqref{pD} contains at most one intersection point with $\sigma^ND$. Additionally, it can contain at most two intersection points with progressions in \eqref{pE}. The latter is due to Lemma \ref{2pt} and the fact that the steps $h_D$ and $h_E$ are linearly independent. As for the rest of points in the progression in \eqref{pD} (which are certainly present with $N>0$ big enough), these are all outside of $\sigma^ND$, and the corresponding monomials in \eqref{ekN(x)} survive under possible reductions. This contradicts to \eqref{ke}.

(c) Consider the case $\operatorname{card}D=2$, that is, $D$ is formed by two distinct points $\binom{i_1}{j_1}$, $\binom{i_2}{j_2}$. Respectively, \eqref{pD} is formed by two progressions whose step is $h_D$ and length is $|a(N)|$. As one observes from \eqref{sigma^N}, the latter value is non-zero and grows for $N>0$ big enough, whose choice is at our hand.

Suppose that these two progressions are disjoint. In this case one can choose either of those to apply the argument in (b) in order to get the desired contradiction. Of course this argument can be also applied to the progressions of \eqref{pE}.

Finally, suppose that the two progressions in \eqref{pD} are not disjoint (at least when their length $|a(N)|$ is big enough). More precisely, one has to assume that $\binom{i_2}{j_2}=\binom{i_1}{j_1}+sh_D$ for some integer $s$. In fact $s=1$, because otherwise the symmetric difference of these two progressions contains $2s$ (at least $4$) points. After discarding at most one intersection point with progressions of \eqref{pE}, we obtain at least $3$ points in \eqref{pD} corresponding to at least $3$ monomials in \eqref{ekN(x)} which survive after all the reductions in \eqref{ekN(x)}. Since \eqref{kNe(x)} contains at most $2$ monomials, we get a contradiction with \eqref{ke}. Thus $s=1$, and a similar argument works also for $E$. \hfill $\blacksquare$

\medskip

Turn back to the {\bf Proof} of Proposition \ref{real_eigens}. At this step we need to diverge from our previous approach based on disregarding the specific non-zero constant multipliers at monomials. Also, we now assume that $\mathsf{k}$ acts via an automorphism $\varphi_{\sigma,\alpha,\beta}$ of the general form. Namely, we have
\begin{equation}\label{k_gen}
\pi(\mathsf{k})(x)=\alpha x^{a(1)}y^{c(1)};\qquad
\pi(\mathsf{k})(y)=\beta x^{b(1)}y^{d(1)}.
\end{equation}
Also, by Lemma \ref{2tf} we have
\begin{align}
\pi(\mathsf{e})(x) &=a_1x^iy^j+a_2x^{i+a(1)-1}y^{j+c(1)},&\label{e(x)}
\\ \pi(\mathsf{e})(y) &=b_1x^{i'}y^{j'}+b_2x^{i'+b(1)}y^{j'+d(1)-1}&\label{e(y)}
\end{align}
for some $a_1,a_2,b_1,b_2\in\mathbb{C}\setminus\{0\}$. Let us compute
\begin{multline}\label{eyx}
\pi(\mathsf{e})(yx)= y\pi(\mathsf{e})(x)+\pi(\mathsf{e})(y)\pi(\mathsf{k})(x)=
\\ =a_1q^ix^iy^{j+1}+a_2q^{i+a(1)-1}x^{i+a(1)-1}y^{j+c(1)+1}+b_1\alpha q^{a(1)j'}x^{i'+a(1)}y^{j'+c(1)}+
\\ +b_2\alpha q^{a(1)(j'+d(1)-1)}x^{i'+a(1)+b(1)}y^{j'+c(1)+d(1)-1},
\end{multline}
\begin{multline}\label{exy}
q\pi(\mathsf{e})(xy)=
qx\pi(\mathsf{e})(y)+q\pi(\mathsf{e})(x)\pi(\mathsf{k})(y)=
\\ =b_1qx^{i'+1}y^{j'}+b_2qx^{i'+b(1)+1}y^{j'+d(1)-1}+a_1\beta q^{b(1)j+1}x^{i+b(1)}y^{j+d(1)}+
\\ +a_2\beta q^{b(1)(j+c(1))+1}x^{i+a(1)+b(1)-1}y^{j+c(1)+d(1)}.
\end{multline}

Equating \eqref{eyx} and \eqref{exy}, we obtain a relation with a sum of $8$ non-zero monomials being equal to zero. Let us consider the picture appearing via the one-to-one correspondence
$\operatorname{const}x^iy^j\mapsto\binom{i}{j}$ between the monomials in
$\mathbb{C}_q[x^{\pm 1},y^{\pm 1}]$ modulo non-zero constant multipliers
and $\mathbb{Z}^2$. Among the $8$ monomials, those $4$ which contain $i$, $j$ in their exponents, correspond to the following $4$ points in $\mathbb{Z}^2$: $u=\binom{i}{j+1}$, $u+h_D$, $u+h_E$, $u+h_D+h_E$. These $4$ points are pairwise distinct, hence the corresponding monomials are linear independent.

In a similar way, the $4$ monomials which contain $i'$, $j'$ in their exponents, correspond to the following $4$ distinct points in $\mathbb{Z}^2$: $v=\binom{i'+1}{j'}$, $v+h_D$, $v+h_E$, $v+h_D+h_E$; hence the corresponding monomials are linear independent. Note that in both cases we have parallelograms whose sides are determined by the same pair of vectors $h_D$, $h_E$. In view of our observations we conclude that the sum of $8$ monomials can be zero only in the case when the two parallelograms coincide. This implies, in particular, that $i=i'+1$, $j+1=j'$.

With the latter conclusion being taken into account, our next step is in equating the coefficients at the corresponding monomials in \eqref{eyx}, \eqref{exy}. This leads to the following system of equations:\\

$\left\{\begin{array}{ll}
a_1q^i &=b_1q
\\ a_2q^{i+a(1)-1} &=-b_1\alpha q^{a(1)j+a(1)}
\\ a_1\beta q^{b(1)j+1} &=-b_2q
\\ a_2\beta q^{b(1)(j+c(1))+1} &=b_2\alpha q^{a(1)(j+d(1))}
\end{array}\right.$\\

\noindent The general solution of this system has the form
$$
a_1=g,\quad a_2=-g\alpha q^{a(1)j},\quad b_1=gq^{i-1},\quad b_2=-g\beta q^{b(1)j},\qquad g\in\mathbb{C}\setminus\{0\}.
$$
This allows one to rewrite \eqref{e(x)} and \eqref{e(y)} as follows:
\begin{align}
\pi(\mathsf{e})(x) &=gx^iy^j-g\alpha q^{a(1)j}x^{i+a(1)-1}y^{j+c(1)},&\label{e(x)'}
\\ \pi(\mathsf{e})(y) &=gq^{i-1}x^{i-1}y^{j+1}-g\beta q^{b(1)j}x^{i+b(1)-1}y^{j+d(1)}.&\label{e(y)'}
\end{align}

Now let us turn back to \eqref{ekN(x)}, \eqref{kNe(x)}, together with the associated pictures in $\mathbb{Z}^2$ \eqref{pE}, \eqref{pD}, both in the case $N=1$, in order to extract more consequences from \eqref{ke} applied to $x$. Looking at \eqref{pD}, we note, in view of Lemma \ref{2tf}, that the entire picture is formed by the two progressions, both of the same length $|a(1)|$ and with the same step $h_D$. These progressions are non-disjoint; more precisely, they are translates of each other in such a way that the initial point of one of those coincides to the second point of another one. This means that the endpoints of the union of these two progressions correspond to the monomials in \eqref{ekN(x)} which survive after possible reductions within the second sum in \eqref{ekN(x)}. To write these two points $u_1,u_2\in\mathbb{Z}^2$ explicitly, one has to substitute $r=0$ and $r=a(1)$ to \eqref{pD}, respectively, and this is clearly independent of the sign of $a(1)$. Namely, we have
$$
u_1=\binom{i+a(1)-1+b(1)c(1)}{j+c(1)d(1)};\qquad u_2=\binom{i+a(1)^2-1+b(1)c(1)}{j+a(1)c(1)+c(1)d(1)}.
$$
It might look all this assumes $a(1)\ne 0$. However this is true only in the above part of the picture. It will become clear in what follows that the case $a(1)=0$, in spite of vanishing the above two progressions, does not break the entire argument.

In a similar way, let us consider the progressions of \eqref{pE}, whose length is $|c(1)|$ and the step is $h_E$, and write explicitly the endpoints $v_1,v_2\in\mathbb{Z}^2$ explicitly substituting $r=0$ and $r=c(1)$ to \eqref{pE} and using the relations between $i,j,i',j'$:
$$
v_1=\binom{i+a(1)-1}{j+c(1)};\qquad v_2=\binom{i+a(1)-1+b(1)c(1)}{j+c(1)d(1)}.
$$
Note that $c(1)\ne 0$, because the matrix $\sigma$ can not be triangular under our assumptions. Again, the corresponding monomials in \eqref{ekN(x)} survive after possible reductions within the first sum in \eqref{ekN(x)}. However, all the four monomials should somehow vanish when equating \eqref{ekN(x)} to \eqref{kNe(x)} due to \eqref{ke}. In view of this, let us compare $u_1$, $u_2$, $v_1$, $v_2$ to the points of $\sigma D$ which correspond to the monomials of \eqref{kNe(x)}. By \eqref{e(x)'}, $\binom{i}{j}\in D$, hence the two points of $\sigma D$ are just $w_1=\sigma\binom{i}{j}$ and $w_2=\sigma\left(\binom{i}{j}+h_D\right)$.

Note that $u_1=v_2$. It follows that $\{u_2,v_1\}=\sigma D$ (obviously, there is no other opportunity). More precisely, since $u_2=u_1+a(1)h_D=u_1+a(1)(\sigma-I)\binom{1}{0}$ and $v_1=v_2-c(1)h_E=v_2-c(1)(\sigma-I)\binom{0}{1}$, one has
$$
u_2-v_1=a(1)(\sigma-I)\binom{1}{0}+c(1)(\sigma-I)\binom{0}{1}=
(\sigma-I)\binom{a(1)}{c(1)}=(\sigma-I)\sigma\binom{1}{0}.
$$
Also one has
$$w_2-w_1=\sigma h_D=\sigma(\sigma-I)\binom{1}{0},$$
and since the matrices $\sigma$, $\sigma-I$ commute, $u_2-v_1=w_2-w_1$. It follows that $v_1=w_1$, $u_2=w_2$. The first of these equalities can be written as
$$\sigma\binom{i}{j}=\binom{i+a(1)-1}{j+c(1)},$$
hence
$$(\sigma-I)\binom{i}{j}=h_D=(\sigma-I)\binom{1}{0},$$
and since $\sigma-I$ is an invertible matrix, we conclude that $i=1$, $j=0$, and this solution is unique. Certainly, the same result can be deduced from $u_2=w_2$. This, together with \eqref{k_gen}, allows one more adjustment of \eqref{e(x)'}, \eqref{e(y)'}:
\begin{alignat}{2}
\pi(\mathsf{e})(x) &=gx-g\alpha x^{a(1)}y^{c(1)} &&=g(\operatorname{id}-\pi(\mathsf{k}))(x),\label{e(x)''}
\\ \pi(\mathsf{e})(y) &=gy-g\beta x^{b(1)}y^{d(1)} &&=g(\operatorname{id}-\pi(\mathsf{k}))(y),\label{e(y)''}
\end{alignat}
for some $g\in\mathbb{C}\setminus\{0\}$.

A routine verification shows that the linear map $\Phi=g(\operatorname{id}-\pi(\mathsf{k}))$ on $\mathbb{C}_q[x^{\pm 1},y^{\pm 1}]$ is subject to the same rule $\Phi(\xi\eta)=\xi\Phi(\eta)+\Phi(\xi)\pi(\mathsf{k})(\eta)$ as $\pi(\mathsf{e})$. This allows one to extend the relations \eqref{e(x)''}, \eqref{e(y)''} from the generators to the entire algebra $\mathbb{C}_q[x^{\pm 1},y^{\pm 1}]$, so that $\pi(\mathsf{e})=g(\operatorname{id}-\pi(\mathsf{k}))$ identically on $\mathbb{C}_q[x^{\pm 1},y^{\pm 1}]$. Since this map clearly commutes with $\pi(\mathsf{k})$, we conclude that \eqref{ke} fails. This contradiction completes the proof of Proposition. \hfill $\blacksquare$

\medskip

\begin{remark}
In view of the results of this Section, the complete list of $U_q\left(\mathfrak{sl}_2\right)$-symmetries on $\mathbb{C}_q[x^{\pm
1},y^{\pm 1}]$ such that $\mathsf{k}$ acts by an automorphism $\varphi_{\sigma,\alpha,\beta}$ with a non-unit matrix $\sigma$, is given by Theorem \ref{sigma=-I}.
\end{remark}

\bigskip

{\center\section{\boldmath$\sigma=I$: the generic case}\label{ws}}

In the case $\sigma=I$, the action of the Cartan element $\mathsf{k}$ is given by multiplication of the generators $x$, $y$ by the {\it weight constants}. We follow \eqref{Auto} in denoting these weight constants by $\alpha$ and $\beta$, respectively. Certainly, monomials form a basis of weight vectors (eigenvectors for $\pi(\mathsf{k})$), and the associated eigenvalues are called {\it weights}.

Let us consider the case when either $\pi(\mathsf{e})$ or $\pi(\mathsf{f})$ is not identically zero. In this setting we describe series of symmetries which we call {\it generic}.

A pair of non-zero complex constants $\alpha$ and $\beta$ which could appear as weight constants for some $U_q(\mathfrak{sl}_2)$-symmetry of $\mathbb{C}_q[x^{\pm 1},y^{\pm 1}]$, can not be arbitrary. In fact, an obvious consequence of \eqref{ke} claims that $\pi(\mathsf{e})$ sends a vector whose weight is $\gamma$ to a vector whose weight is $q^2\gamma$. In particular, $\pi(\mathsf{e})(x)$, if non-zero, is a sum of monomials with (the same) weight $q^2\alpha$. Since the weight of the monomial $ax^iy^j$ (with $a\ne 0$) is $\alpha^i\beta^j$, one has that $\alpha^u\beta^v=q^2$ for some integers $u,v$. Of course similar conclusions can be also derived by applying \eqref{ke}, \eqref{kf} to $x$ and $y$. Under our assumptions, this argument should work at least once.

The following Theorem covers all but a countable family of admissible pairs of weight constants.

\medskip

\begin{theorem}\label{gener}
Let $\alpha,\beta\in\mathbb{C}\setminus\{0\}$ be such that $\alpha^u\beta^v=q^2$ for some $u,v\in\mathbb{Z}$ and $\alpha^m\ne\beta^n$ for non-zero integers $m$, $n$. Then there exists a one-parameter ($a\in\mathbb{C}\setminus\{0\}$) family of $U_q(\mathfrak{sl}_2)$-symmetries of $\mathbb{C}_q[x^{\pm 1},y^{\pm 1}]$:
\begin{flalign}
\pi(\mathsf{k})(x) &=\alpha x &\pi(\mathsf{k})(y)&=\beta y & \label{gener_k}
\\ \pi(\mathsf{e})(x) &=aq^{uv+3}\frac{1-\alpha q^v}{(1-q^2)^2}x^{u+1}y^v &
\pi(\mathsf{e})(y) &=aq^{uv+3}\frac{q^u-\beta}{(1-q^2)^2}x^uy^{v+1} & \label{gener_e}
\\ \pi(\mathsf{f})(x) &=-\frac{(\alpha^{-1}-q^{-v})}{a}x^{-u+1}y^{-v} &
\pi(\mathsf{f})(y) &=-\frac{(\beta^{-1}q^{-u}-1)}{a}x^{-u}y^{-v+1} &
\label{gener_f}
\end{flalign}
There exist no other symmetries with the weight constants $\alpha$ and $\beta$.
\end{theorem}

\begin{lemma}\label{ratios}
Let $\pi$ be a $U_q(\mathfrak{sl}_2)$-symmetry on $\mathbb{C}_q[x^{\pm 1},y^{\pm 1}]$ such that
\begin{align}
\pi(\mathsf{k})(x) &=\alpha x & \pi(\mathsf{k})(y) &=\beta y &\label{kx_y}
\\ \pi(\mathsf{e})(x) &=\sum_{i,j}a_{i,j}x^iy^j &
\pi(\mathsf{e})(y) &=\sum_{i,j}b_{i,j}x^iy^j &\label{ex_y}
\\ \pi(\mathsf{f})(x) &=\sum_{i,j}c_{i,j}x^iy^j &
\pi(\mathsf{f})(y) &=\sum_{i,j}d_{i,j}x^iy^j, &\label{fx_y}
\end{align}
with $\alpha,\beta\in\mathbb{C}\setminus\{0\}$, $a_{i,j},b_{i,j},c_{i,j},d_{i,j}\in\mathbb{C}$, and the above sums being finite. Then
\begin{align}
a_{i+1,j}\left(q^i-\beta\right) &=b_{i,j+1}\left(1-\alpha q^j\right),\label{ab}
\\ c_{i+1,j}\left(1-\beta^{-1}q^i\right) &=d_{i,j+1}\left(q^j-\alpha^{-1}\right).\label{cd}
\end{align}
\end{lemma}

{\bf Proof.} This is a consequence of \eqref{qpr}. We have
$$
\pi(\mathsf{e})(yx)=
y\pi(\mathsf{e})(x)+\pi(\mathsf{e})(y)\pi(\mathsf{k})(x)=
\sum_{i,j}a_{i,j}q^ix^iy^{j+1}+\sum_{i,j}b_{i,j}\alpha q^jx^{i+1}y^j,
$$
$$
q\pi(\mathsf{e})(xy)=
qx\pi(\mathsf{e})(y)+q\pi(\mathsf{e})(x)\pi(\mathsf{k})(y)=
\sum_{i,j}b_{i,j}qx^{i+1}y^j+\sum_{i,j}a_{i,j}\beta qx^iy^{j+1}.
$$
With this, we project the relation $\pi(\mathsf{e})(yx)=q\pi(\mathsf{e})(xy)$ to the one dimensional subspace $\mathbb{C}x^{i+1}y^{j+1}$ parallel to the linear span of all other monomials to obtain \eqref{ab}. The proof of \eqref{cd} is similar. \hfill $\blacksquare$

\medskip

{\bf Proof} of Theorem \ref{gener}. We are about to apply the following relations valid under the assumption $\sigma=I$ of the present Section:
\begin{align*}
\pi(\mathsf{e})(x^p) &=\sum_{i,j}a_{i,j}\frac{\alpha^pq^{jp}-1}{\alpha q^j-1}x^{p-1+i}y^j, &
\\ \pi(\mathsf{e})(y^p)
&=\sum_{i,j}b_{i,j}\frac{\beta^p-q^{ip}}{\beta-q^i}x^iy^{p-1+j}, &
\\ \pi(\mathsf{f})(x^p)
&=\sum_{i,j}c_{i,j}\frac{\alpha^{-p}-q^{jp}}{\alpha^{-1}-q^j}x^{p-1+i}y^j, &
\\ \pi(\mathsf{f})(y^p)
&=\sum_{i,j}d_{i,j}\frac{\beta^{-p}q^{ip}-1}{\beta^{-1}q^i-1}x^iy^{p-1+j}, &
\end{align*}
where $p\in\mathbb{Z}$, $a_{i,j},b_{i,j},c_{i,j},d_{i,j}\in\mathbb{C}$ are as in \eqref{ex_y}, \eqref{fx_y}. These relations are due to a straightforward induction argument.

A direct computation that applies the above relations, allows one to verify that the extended action \eqref{gener_k} -- \eqref{gener_f} from the generators to the entire algebras $U_q(\mathfrak{sl}_2)$ and $\mathbb{C}_q[x^{\pm 1},y^{\pm 1}]$ passes through all the relations \eqref{qpr}, \eqref{Lpr}, \eqref{kk1} -- \eqref{effe}. Hence \eqref{gener_k} -- \eqref{gener_f} determine a well-defined $U_q(\mathfrak{sl}_2)$-symmetry on $\mathbb{C}_q[x^{\pm 1},y^{\pm 1}]$.

Let us prove that there are no other symmetries. Observe first that the assumptions of the Theorem on the weight constants $\alpha$ and $\beta$ imply that the pair of integers $u,v$ with $\alpha^u\beta^v=q^2$ is unique. Therefore the monomials modulo a (non-zero) constant multiplier are in one-to-one correspondence with their weights. Since, in view of \eqref{ke}, \eqref{kf}, $\pi(\mathsf{e})$ and $\pi(\mathsf{f})$ `multiply' the weight of a weight vector by $q^2$ and $q^{-2}$, respectively, one deduces that $\pi(\mathsf{e})(x)$, $\pi(\mathsf{e})(y)$, $\pi(\mathsf{f})(x)$, $\pi(\mathsf{f})(y)$ are monomials which, up to constant multipliers, should be just as in \eqref{gener_e}, \eqref{gener_f}.

Observe that, under the assumptions of the Theorem on weight constants $\alpha$ and $\beta$, no differences in \eqref{gener_e}, \eqref{gener_f} could be zero. Thus it follows from Lemma \ref{ratios} that the ratio of coefficients at the monomials $\pi(\mathsf{e})(x)$ and $\pi(\mathsf{e})(y)$ should be just as in \eqref{gener_e}. Of course, a similar claim is also true for \eqref{gener_f}.

It remains to establish the ratio of coefficients in \eqref{gener_e} and those in \eqref{gener_f}. This is done via applying \eqref{effe} to $x$ and $y$. The computation in question, which is left to the reader, in both cases, leads to the same result reflected in \eqref{gener_e}, \eqref{gener_f}. \hfill $\blacksquare$

\begin{remark}
Theorem \ref{gener} describes an uncountable family of isomorphism classes of $U_q(\mathfrak{sl}_2)$-symmetries of $\mathbb{C}_q[x^{\pm 1},y^{\pm 1}]$.

In fact, one clearly has an uncountable family of admissible (i.e., those subject to the assumptions of Theorem \ref{gener}) pairs of weight constants $\alpha,\beta$. On the other hand, the action of the group of automorphisms of $\mathbb{C}_q[x^{\pm 1},y^{\pm 1}]$ (which is the semidirect product of its subgroups $SL(2,\mathbb{Z})$ and $(\mathbb{C}^*)^2$) on the space of parameters of generic symmetries is such that the action of normal subgroup $(\mathbb{C}^*)^2$ remains intact every pair of weight constants $(\alpha,\beta)$. It follows that (the projection of) each orbit of the automorphism group on the space of admissible pairs of weight constants $\alpha,\beta$ is only countable.
\end{remark}


\bigskip

\begin{thebibliography}{99}

\bibitem{abe} {\it E. Abe}, Hopf Algebras. Cambridge Univ. Press,
    Cambridge, 1980.

\bibitem{AD} {\it J. Alev and F. Dumas}, Rigidit\'e des plongements des quotients primitifs minimaux de $U_q(\mathfrak{sl}(2))$ dans l'Alg\`ebre quantique de Weyl-Hayashi. -- Nagoya Math. J. {\bf 143} (1996), 119 -- 146.

\bibitem{DHL} {\it S. Duplij, Y. Hong, and F. Li}, Analysis of
    $U_q(\mathfrak{sl}_{m+1})$-symmetries on quantum $n$-spaces. --
    arXiv:1305.6582 [math.QA], 37 p. (To appear in J. Lie Theory)

\bibitem{DS} {\it S. Duplij and S. Sinel'shchikov}, Classification of
    $U_q(\mathfrak{sl}_2)$-module algebra structures on the quantum plane,
    -- J. Math. Phys., Anal., and Geom. {\bf 6} (2010), No 4, 406 -- 430.

\bibitem{kassel} {\it C. Kassel}, Quantum Groups. Springer-Verlag, New
    York, 1995, 531 p.p.

\bibitem{KPS} {\it E. Kirkman, C. Procesi, L. Small}, A q-Analog for the
    Virasoro Algebra. -- {\it Comm. Algebra} {\bf 22} (10), 3755 -- 3774.

\bibitem{PLCCN} {\it Park Hong Goo, Lee Jeongsig, Choi Seul Hee, Chen
    XueQing, and Nam Ki-Bong}, Automorphism groups of some algebras. --
    Science in China Series A: Mathematics {\bf 52} (2009), No 2, 323 --
    328.

\bibitem{sweedler} {\it M. E. Sweedler}, Hopf Algebras. Benjamin, New York,    1969.

\end{thebibliography}
\end{document}